# THE COMPOSITE ABSOLUTE PENALTIES FAMILY FOR GROUPED AND HIERARCHICAL VARIABLE SELECTION[1]

By Peng Zhao, Guilherme Rocha and Bin Yu[2]

*University of California, Berkeley*

Extracting useful information from high-dimensional data is an important focus of today's statistical research and practice. Penalized loss function minimization has been shown to be effective for this task both theoretically and empirically. With the virtues of both regularization and sparsity, the $L_1$-penalized squared error minimization method Lasso has been popular in regression models and beyond.

In this paper, we combine different norms including $L_1$ to form an intelligent penalty in order to add side information to the fitting of a regression or classification model to obtain reasonable estimates. Specifically, we introduce the Composite Absolute Penalties (CAP) family, which allows given grouping and hierarchical relationships between the predictors to be expressed. CAP penalties are built by defining groups and combining the properties of norm penalties at the across-group and within-group levels. Grouped selection occurs for nonoverlapping groups. Hierarchical variable selection is reached by defining groups with particular overlapping patterns. We propose using the BLASSO and cross-validation to compute CAP estimates in general. For a subfamily of CAP estimates involving only the $L_1$ and $L_\infty$ norms, we introduce the iCAP algorithm to trace the entire regularization path for the grouped selection problem. Within this subfamily, unbiased estimates of the degrees of freedom (df) are derived so that the regularization parameter is selected without cross-validation. CAP is shown to improve on the predictive performance of the LASSO in a series of simulated experiments, including cases with $p \gg n$ and possibly mis-specified groupings. When the complexity of a model is properly calculated, iCAP is seen to be parsimonious in the experiments.

Received March 2007; revised October 2007.

[1]Supported in part by NSF Grants DMS-06-05165, DMS-03-036508, DMS-04-26227, ARO Grant W911NF-05-1-0104 and NSFC Grant 60628102.

[2]Supported by a Miller Research Professorship in Spring 2004 from the Miller Institute at UC Berkeley and a Guggenheim Fellowship in 2006.

**1. Introduction.** Information technology advances are bringing the possibility of new and exciting discoveries in various scientific fields. At the same time, they pose challenges for the practice of statistics because current data sets often contain a large number of variables compared to the number of observations. A paramount example is micro-array data, where thousands or more gene expressions are collected while the number of samples remains around a few hundreds (e.g., [11]).

For regression and classification, parameter estimates are often defined as the minimizer of an empirical loss function $L$ and too responsive to noise when the number of observations $n$ is small with respect to the model dimensionality $p$. Regularization procedures impose constraints represented by a penalization function $T$. The regularized estimates are given by:

$$\hat{\beta}(\lambda) = \arg\min_{\beta}[L(Z,\beta) + \lambda \cdot T(\beta)],$$

where $\lambda$ controls the amount of regularization and $Z = (X, Y)$ is the (random) observed data: $X$ is the $n \times p$ design matrix of predictors, and $Y$ a $n$-dimensional vector of response variables ($Y$ being continuous for regression and discrete for classification). We restrict attention to linear models:

$$(1) \qquad L(Z, \beta) = L\left(Y, \sum_{j=1}^{p} \beta_j X_j\right),$$

where $X_j$ denotes the observed values for the $j$th predictor ($j$th column of $X$). Thus, setting $\beta_j = 0$ corresponds to excluding $X_j$ from the model.

Selection of variables is a popular way of performing regularization. It stabilizes the parameter estimates while leading to interpretable models. Early variable selection methods such as Akaike's *AIC* [1], Mallow's $C_p$ [15] and Schwartz's *BIC* [19] are based on penalizing the dimensionality of the model. Penalizing estimates by their Euclidean norm (ridge regression [12]) is commonly used among statisticians. Bridge estimates [9] use the $L_\gamma$-norm on the $\beta$ parameter defined as $\|\beta\|_\gamma = (\sum_{j=1}^p |\beta_j|^\gamma)^{1/\gamma}$ as a penalization: they were first considered as a unifying framework in which to understand ridge regression and variable selection (the dimensionality of the model being interpreted as the "$L_0$-norm" of the coefficients). More recently, the nonnegative garrote [3], wavelet shrinkage [5], basis pursuit [4] and the LASSO [22] have exploited the convexity [2] of the $L_1$-norm as a more computationally tractable means for selecting variables.

For severely ill-posed estimation problems, sparsity alone may not be sufficient to obtain stable estimates [26]. Group or hierarchical information can be a source of further regularization constraints. Sources of such information vary according to the problem at hand. A grouping of the predictors may arise naturally: for categorical variables, the dummy variables used to represent different levels define a natural grouping [23]. Alternatively, a natural hierarchy may exist: an interaction term is usually only included in a model after its corresponding main effects. More broadly, in applied work, expert knowledge is a potential source of grouping and hierarchical information. After completing this work, we have also learned that grouping of parameters also occur naturally when fitting various regressions simultaneously [16].

In this paper, we introduce the Composite Absolute Penalties (CAP) family of penalties. CAP penalties are highly customizable and build upon $L_\gamma$ penalties to express both grouped and hierarchical selection. The overlapping patterns of the groups and the norms applied to the groups of coefficients used to build a CAP penalty determine the properties of the associated estimate. The CAP penalty formation assumes a known grouping or hierarchical structure on the predictors. For group selection, clustering techniques can be used to define groups of predictors before the CAP penalty is applied. In our simulations, this two-step approach has resulted in CAP estimates that were robust to misspecified groups.

Zou and Hastie [26] (the Elastic Net), Kim, Kim and Kim [14] (Blockwise Sparse Regression) and Yuan and Lin [23] (GLASSO) have previously explored combinations of the $L_1$-norm and $L_2$-norm penalties to achieve more structured estimates. The CAP family extends these ideas in two directions: first, it allows different norms to be combined and; second, different overlapping of the groups are allowed to be used. These extensions result both in computational and modeling gains. By allowing norms other than $L_1$ and $L_2$ to be used, the CAP family allows computationally convenient penalties to be constructed from the $L_1$ and $L_\infty$ norms. By letting the groups overlap, CAP penalties can be constructed to represent a hierarchy among the predictors. In Section 2.2.2, we detail how the groups should overlap for a given hierarchy to be represented.

CAP penalties built from the $L_1$ and $L_\infty$-norms are computationally convenient, as their regularization paths are piecewise linear for piecewise quadratic loss functions [18]. We call such group of penalties the iCAP family ("i" standing for the infinity norm). We extend the homotopy/LARS-LASSO algorithm [8, 17] and design fast algorithms for iCAP penalties in the cases of nonoverlapping group selection (the iCAP algorithm) and tree-hierarchical selection (the hierarchical iCAP algorithm: hiCAP). A Matlab implementation of these algorithms is available from: http://www.stat.berkeley.edu/twiki/Research/YuGroup/Software.

For iCAP penalties, used with the $L_2$-loss, unbiased estimates of the degrees of freedom of models along the path can be obtained by extending the results in Zou, Hastie and Tibshirani [27]. It is then possible to employ information theory criteria to pick an estimate from the regularization path and thus avoiding the use of cross validation. Models picked from the iCAP path using Sugiura's [21] $AIC_c$ and the degrees of freedom estimates had predictive performance comparable to cross-validated models even when $n \ll p$.

The computational advantage of CAP penalties is preserved in a broader setting. We prove that CAP penalties is convex whenever all norms used in its construction are convex. Based on this, we propose using the BLASSO algorithm [24] to compute the CAP regularization path and cross-validation to select the amount of regularization.

Our experimental results show that the inclusion of group and hierarchical information substantially enhance the predictive performance of the penalized estimates in comparison to the LASSO. This improvement was preserved even when empirically determined partitions of the set of predictors was severely mis-specified and was observed for different settings of the norms used to build CAP. While the CAP estimates are not sparser than LASSO estimates in the number of variables sense, they result in more parsimonious use of degrees of freedom and more stable estimates [7].

The remainder of this paper is organized as follows. In Section 2, for a given grouping or hierarchical structure, we define CAP penalties, relate them to the properties of $\gamma$-norm penalized estimates and detail how to build CAP penalties from given group and hierarchical information. Section 3 proves the convexity of the CAP penalties and describes the computation of CAP estimates. We propose algorithms for tracing the CAP regularization path and methods for selecting the regularization parameter $\lambda$. Section 4 gives experimental results based on simulations of CAP regression with the $L_2$-loss and explores a data-driven group formation procedure showing that CAP estimates enjoy some robustness relative to possibly mis-specified groupings. Section 5 concludes the paper with a brief discussion.

**2. The Composite Absolute Penalty (CAP) family.** We first give a review of the properties of $L_\gamma$-norm penalized estimates. Then we show how CAP penalties exploit them to reach grouped and hierarchical selection. For overlapping groups, our focus will be on how to overlap groups so hierarchical selection is achieved.

2.1. *Preliminaries: properties of bridge regressions.* We consider an extended version of the bridge regression [9] where a general loss function replaces the $L_2$-loss. The bridge regularized coefficients are given by

(2) $$\hat{\beta}_\gamma(\lambda) = \arg\min_\beta [L(Z, \beta) + \lambda \cdot T(\beta)]$$

$$\text{with } T(\beta) = \|\beta\|_\gamma^\gamma = \sum_{j=1}^p |\beta_j|^\gamma.$$

The properties of bridge estimates path vary considerably according to the value chosen for $\gamma$. The estimates tend to fall in regions of high "curvature" of the penalty contour plot: for $0 \leq \gamma \leq 1$, some estimated coefficients are set to zero; for $1 < \gamma < 2$, estimated coefficients lying in directions closer to the axis are favored; for $\gamma = 2$, the estimates are not encouraged to lie in any particular direction; and finally for $2 < \gamma \leq \infty$, they tend to concentrate along the diagonals. Figure 1 illustrates the different behavior of bridge regression estimates for different values of $\gamma$ for the diabetes data used in Efron et al. [8].

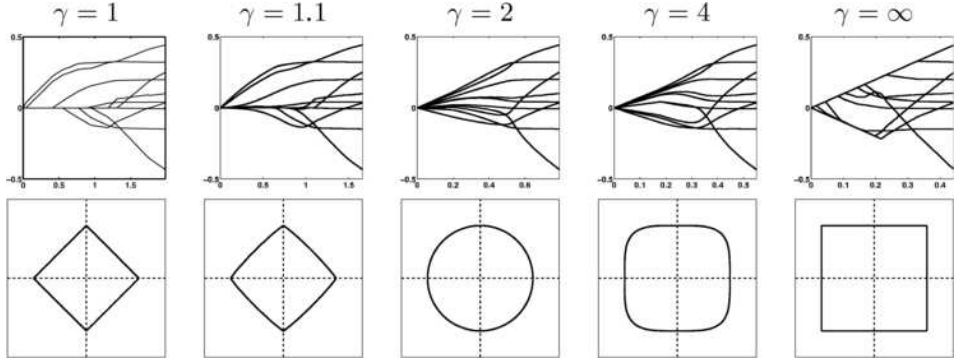

FIG. 1. **Regularization paths of bridge regressions.** *Upper Panel: Regularization paths for different bridge parameters for the diabetes data. From left to right: Lasso ($\gamma = 1$), near-Lasso ($\gamma = 1.1$), Ridge ($\gamma = 2$), over-Ridge ($\gamma = 4$), max($\gamma = \infty$). The horizontal and vertical axis contain the $\gamma$-norm of the normalized coefficients and the normalized coefficients respectively. Lower Panel: Contour plots $\|(\beta_1, \beta_2)\|_\gamma = 1$ for the corresponding penalties.*

CAP penalties exploit the distinct behaviors of the bridge estimates according to whether $\gamma = 1$ or $\gamma > 1$. For convex $L$ and $\gamma \geq 1$, the bridge estimates are fully characterized by the Karush–Kuhn–Tucker (KKT) conditions (see [2]):

$$\text{(3)} \quad \frac{\partial L}{\partial \beta_j} = -\lambda \frac{\partial \|\beta\|_\gamma}{\partial \beta_j} = -\lambda \cdot \text{sign}(\beta_j) \frac{|\beta_j|^{\gamma-1}}{\|\beta\|_\gamma^{\gamma-1}} \qquad \text{for } j \text{ such that } \beta_j \neq 0;$$

$$\text{(4)} \quad \left|\frac{\partial L}{\partial \beta_j}\right| \leq \lambda \left|\frac{\partial \|\beta\|_\gamma}{\partial \beta_j}\right| = \lambda \frac{|\beta_j|^{\gamma-1}}{\|\beta\|_\gamma^{\gamma-1}} \qquad \text{for } j \text{ such that } \beta_j = 0.$$

Hence, for $1 < \gamma \leq \infty$, the estimate $\hat{\beta}_j$ equals zero if and only if $\frac{\partial L(Y_i, X_i, \hat{\beta})}{\partial \beta_j}|_{\beta_j=0} = 0$. This condition is satisfied with probability zero when the distribution of $Z_i = (X_i, Y_i)$ is continuous and L is strictly convex. Therefore, $1 < \gamma \leq \infty$ implies that all variables are almost surely included in the bridge estimate. When $\gamma = 1$, however, the right-hand side of (4) becomes a constant set by $\lambda$ and thus variables that fail to infinitesimally reduce the loss by a certain threshold are kept at zero. In what follows, we show how these distinctive behaviors result in group and hierarchical selections.

2.2. *CAP penalties.* We start this subsection by defining the CAP family in its most general form. Then for a given group or hierarchical structure, we specialize the CAP penalty for grouped and hierarchical selection. Unless otherwise stated, we assume that each predictor $X_j$ in what follows is normalized to have mean zero and variance one.

Let
$$\mathcal{I} = \{1, \ldots, p\}$$
contain all indices of the predictors. Given $K$ subsets of indices,
$$\mathcal{G}_k \subset \mathcal{I}.$$
The group formation varies according to the given grouping or hierarchical structure that we want to express through the CAP penalty. Details are presented later in this section. We let a given grouping be denoted by
$$\mathcal{G} = (\mathcal{G}_1, \ldots, \mathcal{G}_K).$$

Moreover, a vector of norm parameters $\gamma = (\gamma_0, \gamma_1, \ldots, \gamma_K) \in \mathbb{R}_+^{K+1}$ must be defined. We let $\gamma_k \equiv c$ denote the case $\gamma_k = c, \forall k \geq 1$.

Call the $L_{\gamma_0}$-norm the overall norm and $L_{\gamma_k}$-norm the $k$th group norm and define

(5)
$$\begin{aligned} \beta_{\mathcal{G}_k} &= (\beta_j)_{j \in \mathcal{G}_k}, \\ N_k &= \|\beta_{\mathcal{G}_k}\|_{\gamma_k} \quad \text{and} \\ \mathbf{N} &= (N_1, \ldots, N_K) \qquad \text{for } k = 1, \ldots, K. \end{aligned}$$

The CAP penalty for grouping $\mathcal{G}$ and norms $\gamma$ is given by

(6)
$$T_{\mathcal{G},\gamma}(\beta) = \|\mathbf{N}\|_{\gamma_0}^{\gamma_0} = \left[\sum_k |N_k|^{\gamma_0}\right].$$

The corresponding CAP estimate for the regularization parameter $\lambda$ is

(7)
$$\hat{\beta}_{\mathcal{G},\gamma}(\lambda) = \arg\min_{\beta}[L(Z,\beta) + \lambda \cdot T_{\mathcal{G},\gamma}(\beta)].$$

In its full generality, the CAP penalties defined above can be used to induce a wide array of different structures in the coefficients: $\gamma_0$ determines how groups relate to one another while $\gamma_k$ dictates the relationship of the coefficients within group $k$. The general principle follows from the distinctive behavior of bridge estimates for $\gamma > 1$ and $\gamma = 1$ as discussed above. Hence, for $\gamma_0 = 1$ and $\gamma_k > 1$ for all $k$, the variables in each group are selected as a block [14, 23]. Nonoverlapping groups yield grouped selection, while suitably constructed overlapping patterns can be used to achieve hierarchical selection. More general overlapping patterns and norm choices are possible, but we defer their study for future research as they are not needed for our goal of grouped and hierarchical selection.

GROUPED AND HIERARCHICAL SELECTION 7

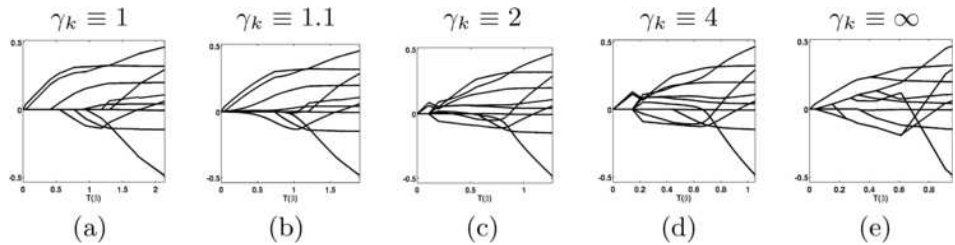

FIG. 2. **Effect of group-norm on regularization path.** *In this figure, we show the regularization path for CAP penalties with different group norms applied to the diabetes data in Efron et al. [8]. The predictors were split into three groups: the first group contains age and sex; the second, body mass index and blood pressure; and the third, blood serum measurements. From left to right, we see:* (a) *Lasso* ($\gamma_0 = 1, \gamma_k \equiv 1$); (b) *CAP*(1.1), ($\gamma_0 = 1, \gamma_k \equiv 1.1$); (c) *GLasso* ($\gamma_0 = 1, \gamma_k \equiv 2$), (d) *CAP*(4), ($\gamma_0 = 1, \gamma_k \equiv 4$); (e) *iCAP*($\gamma_0 = 1, \gamma_k \equiv \infty$).

2.2.1. *Grouped selection: the nonoverlapping groups case.* When the goal of the penalization is to select or exclude nonoverlapping groups of variables simultaneously and the groups are known, we form nonoverlapping groups $\mathcal{G}_k, k = 1, \ldots, K$ to reflect this information. That is, all variables to be added or deleted concurrently should be collected in one group $\mathcal{G}_k \in \mathcal{G}$.

Given the grouping $\mathcal{G}$, the CAP penalization can be interpreted as mimicking the behavior of bridge penalties on two different levels: an across-group and a within-group level. On the across-group level, the group norms $N_k$ behave as if they were penalized by a $L_{\gamma_0}$-norm. On the within-group level, the $\gamma_k$ norm then determines how the coefficients $\beta_{\mathcal{G}_k}$ relate to each other. A formal result establishing this is established in a Bayesian interpretation of CAP penalties for nonoverlapping groups presented in details in the technical report version of this paper [25]. For $\gamma_0 = 1$, sparsity in the **N** vector of group norms is promoted. The $\gamma_k > 1$ parameters then determine how close together the size of the coefficients within a selected group are kept. Thus, Yuan and Lin's [23] corresponds to the LASSO on the across-group level and the rotational invariant ridge penalization on the within-group level. Figure 2 illustrates this fact for the diabetes data from Efron et al. [8].

By allowing group norms other than $L_2$ to be applied to the coefficients in a group, CAP can lead to computational savings by setting $\gamma_0 = 1$ and $\gamma_k \equiv \infty$. In Section 3 below, we present computationally efficient algorithm and model selection criterion for such CAP penalties.

In Section 4, simulation experiments provide compelling evidence that the addition of the group structure can greatly enhance the predictive performance of an estimated model.

*A note on normalization.* Assuming the predictors are normalized, it may be desirable to account for the size of different groups when building the penalties. For $\gamma_0 = 1$, $\gamma_k \equiv \bar{\gamma}$ and letting $\bar{\gamma}^* = \frac{\bar{\gamma}}{\bar{\gamma}-1}$, the decision on whether group $\mathcal{G}_k$ is in the model for a fixed $\lambda$ can be shown to depend upon $\|\nabla_{\beta_{\mathcal{G}_k}} L(Z, \hat{\beta}(\lambda))\|_{\bar{\gamma}^*}$. Thus, larger groups are more likely to be included in the model purely due to their size. We propose that group normalization is achieved by dividing the variance normalized predictors in $\mathcal{G}_k$ by $q_k^{1/\bar{\gamma}^*}$. Following Yuan and Lin [23], such correction causes two hypothetical groups $\mathcal{G}_k$ and $\mathcal{G}_{k'}$ having $\|\nabla_{\beta_j} L(Z, \hat{\beta}(\lambda))\| = c$ for all $j \in \mathcal{G}_k \cup \mathcal{G}_{k'}$ to be included in the CAP estimate simultaneously. Note that for $\bar{\gamma} = 1$ (LASSO), $\bar{\gamma}^* = \infty$ and the group sizes are ignored as in this setting the group structure is lost.

In an extended technical report version of this paper [25], we perform experiments suggesting that the additional normalization by group size does not affect the selection results greatly. This aspect of grouped CAP penalties is explained in further detail there.

2.2.2. *Hierarchical selection: the nested groups case.* We now show how to define CAP penalties to achieve hierarchical selection. We start from an example used to illustrate the principle behind hierarchical CAP penalties, and then prove a result concerning CAP penalties with overlapping groups. Finally, we show how to use our result to build a penalty that induces a hierarchy starting from its representation as a directed graph.

Consider a simple case involving two predictors $X_1$ and $X_2$. Suppose we want $X_1$ to be included in the model before $X_2$. A directed graph can be used to represent this hierarchy as shown in panel (a) in Figure 4. This hierarchy can be induced by defining the overlapping groups $\mathcal{G}_1 = \{1, 2\}$ and $\mathcal{G}_2 = \{2\}$ with $\gamma_0 = 1$, $\gamma_m > 1$ for $m = 1, 2$. That results in the penalty

$$(8) \qquad T(\beta) = \|(\beta_1, \beta_2)\|_{\gamma_1} + \|(\beta_2)\|_{\gamma_2}.$$

The contour plots of this penalty function are shown in Figure 3 for different values of $\gamma = \gamma_1$. As $\mathcal{G}_2$ contains only one variable, these contours are the same regardless of the value chosen for $\gamma_2$ ($\|\beta_2\|_{\gamma_2} = |\beta_2|$ for any $\gamma_2$). The breakpoints along the $\beta_2 = 0$ axis in panels (b) through (d) show that solutions with $\beta_1 \neq 0$ and $\beta_2 = 0$ tend to be encouraged by this penalty when $\gamma_1 > 1$. Setting $\gamma_1 = 1$, however, causes breakpoints to appear along the $\beta_1 = 0$ axis as shown in panel (a), hinting that $\gamma_1 > 1$ is needed for the hierarchical structure to be preserved.

In what refers to the definition of the groups, two things were important for the hierarchy to arise from penalty (8): first, $\beta_2$ was in every group $\beta_1$ was; second, there was one group in which $\beta_2$ was penalized without $\beta_1$ being penalized. As we will see below, having $\beta_2$ in every group where $\beta_1$ is ensures that, once $\beta_2$ deviates from zero, the infinitesimal penalty of $\beta_1$ becomes zero. In addition, letting $\beta_2$ be on a group of its own makes it possible for $\beta_1$ to deviate from zero, while $\beta_2$ is kept at zero.

*The general principle behind the example.* The example above suggests that, to construct more general hierarchies, the key is to set $\gamma_0 = 1$, $\gamma_k > 1$ for all $k$. Given such a $\gamma$, a penalty can cause a set of indices $\mathcal{I}_1$ to be added before the set $\mathcal{I}_2$, by defining groups $\mathcal{G}_1 = \mathcal{I}_2$ and $\mathcal{G}_2 = \mathcal{I}_1 \cup \mathcal{I}_2$. Our next result extends this simple case to more interesting hierarchical structures.

THEOREM 1. *Let $\mathcal{I}_1, \mathcal{I}_2 \subset \{1, \ldots, p\}$ be two subsets of indices. Suppose:*

- $\gamma_0 = 1$ and $\gamma_k > 1, \forall k = 1, \ldots, K$.
- $\mathcal{I}_1 \subset \mathcal{G}_k \Rightarrow \mathcal{I}_2 \subset \mathcal{G}_k$ for all $k$ and
- $\exists k^*$ such that $\mathcal{I}_2 \subset \mathcal{G}_{k^*}$ and $\mathcal{I}_1 \not\subset \mathcal{G}_{k^*}$.

*Then, $\frac{\partial}{\partial \beta_{\mathcal{I}_1}} T(\beta) = 0$ whenever $\beta_{\mathcal{I}_2} \neq 0$ and $\beta_{\mathcal{I}_1} = 0$.*

A proof is given in the Appendix A. Assuming the set $\{Z \in \mathbb{R}^{n \times (p+1)} : \frac{\partial}{\partial \beta_{\mathcal{I}_1}} L(Z, \beta)|_{\beta_{\mathcal{I}_2} = \hat{\beta}_{\mathcal{I}_2}, \beta_{\mathcal{I}_1} = 0} = 0\}$ to have zero probability, Theorem 1 states that once the variables in $\mathcal{I}_2$ are added to the model, infinitesimal movements of the coefficients of variables in $\mathcal{I}_1$ are not penalized and hence $\beta_{\mathcal{I}_1}$ will almost surely deviate from zero.

*Defining groups for hierarchical selection.* Using Theorem 1, a grouping for a more complex hierarchical structure can be constructed from its representation as a directed graph. Let each node correspond to a group of variables $\mathcal{G}_k$ and set its descendants to be the groups that should only be added to the model after $\mathcal{G}_k$. The graph representing the hierarchy in the simple case with two predictors above is shown in the panel (a) of Figure 4. Based on this representation and Theorem 1, CAP penalties enforcing the given hierarchy can be obtained by setting

$$
(9) \qquad T(\beta) = \sum_{m=1}^{\text{nodes}} \alpha_m \cdot \|(\beta_{\mathcal{G}_m}, \beta_{\text{all descendants of } \mathcal{G}_m})\|_{\gamma_m}
$$

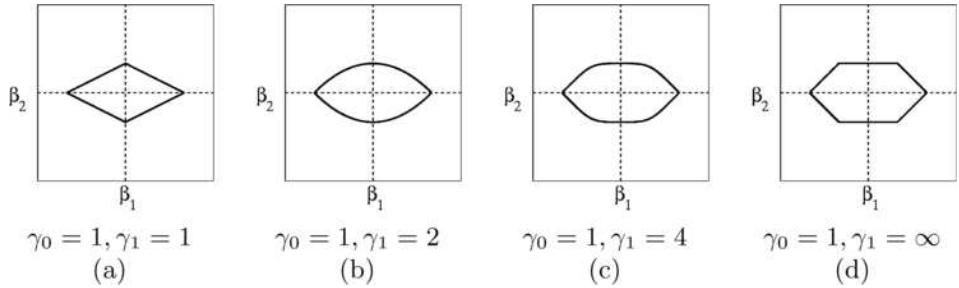

FIG. 3. *Contour plots for the penalty in (8).*

with $\alpha_m > 0$ for all $m$. The factor $\alpha_m$ can be used to correct for the effect of a coefficient being present in too many groups, a concern brought to our attention by one of the referees. In this paper, we keep $\alpha_m = 1$ for all $m$ throughout. We return to this issue in our experimental section below.

For a more concrete example, consider a regression model involving $d$ predictors $x_1, \ldots, x_d$ and all its second-order interactions. Suppose that an interaction term $x_i x_j$ is to be added only after the corresponding main effects $x_i$ and $x_j$. The hierarchy graph is formed by adding an arrow from each main effect to each of its interaction terms. Figure 4 shows the hierarchy graph for $d = 4$. Figure 5 shows sample paths for $d = 4$ with penalties based on (9) and using this hierarchy graph and having different settings for $\gamma_k$. These sample paths were obtained by setting $\beta_1 = 20$, $\beta_2 = 10$, $\beta_3 = 5$, $\beta_{1,2} = 15$ and $\beta_{2,3} = 7$. All remaining coefficients in the model are set to zero. Each predictor has standard normal distribution and the signal to noise ratio is set to 2. Because of the large effect of some interaction terms they are added to the model before their respective main effects when the LASSO is used. However, setting $\gamma_k$ to be slightly larger than 1 is already enough to cause the hierarchical structure to be satisfied. We develop this example further in one of the simulation studies in Section 4.

**3. Computing CAP estimates.** The proper value of the regularization parameter $\lambda$ to use with CAP penalties is rarely known in advance. Two ingredients are then needed to implement CAP in practice: efficient ways of computing estimates for different values of $\lambda$ and a criterion for choosing an appropriate $\lambda$. This section proposes methods for completing these tasks.

3.1. *Tracing the CAP regularization path.* Convexity is a key property for solving optimization problems such as the one defining CAP estimates (7). When the objective function is convex, a point satisfying the Karush–Kuhn–Tucker (KKT) conditions is necessarily a global minimum (see [2]). As the algorithms we present below rely on tracing solutions to the KKT conditions for different values of $\lambda$, we now present sufficient conditions for convexity of the CAP program (7). A proof is given in Appendix A.

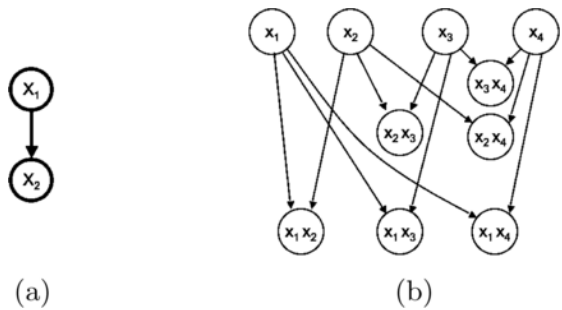

FIG. 4. **Directed graphs representing hierarchies:** (a) *The hierarchy of the simple example: $X_1$ must precede $X_2$;* (b) *Hierarchy for a main and interaction effects model with four variables. The "root nodes" correspond to the main effects and must be added to the model before its children. Each main effect has all the interactions in which it takes part as its children. Each second order interaction effect has two parents.*

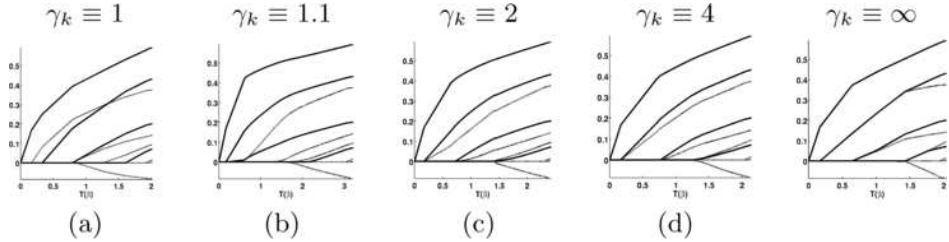

FIG. 5. **A sample regularization path for the simple ANOVA example with four variables.** *In the LASSO path, an interaction (dotted lines) is allowed to enter the model when its corresponding main effect (solid lines) is not in the model. When the group norm $\gamma_k$ is greater than one, the hierarchy is respected. From left to right:* (a) *Lasso* ($\gamma_0 = 1, \gamma_k \equiv 1$); (b) *CAP*(1.1), ($\gamma_0 = 1, \gamma_k \equiv 4$), (c) *GLasso* ($\gamma_0 = 1, \gamma_k \equiv 2$), (d) *CAP*(4), ($\gamma_0 = 1, \gamma_k \equiv 4$), (e) *iCAP*, ($\gamma_0 = 1, \gamma_k \equiv \infty$).

THEOREM 2. *If $\gamma_i \geq 1, \forall i = 0, \ldots, K$, then $T(\beta)$ in (6) is convex. If, in addition, the loss function $L$ is convex in $\beta$ the objective function of the CAP optimization problem in (7) is convex.*

We now detail algorithms for computing the CAP regularization path. The BLasso algorithm is used to deal with general convex loss functions and CAP penalties. Under the $L_2$-loss with $\gamma_0 = 1$ and $\gamma_k \equiv \infty$, we introduce the iCAP ($\infty$-CAP) and the hiCAP (hierarchical-$\infty$-CAP) algorithms to trace the path for group and tree-structured hierarchical selection, respectively.

3.1.1. *The BLasso algorithm.* The BLasso algorithm [24] can be used to approximate the regularization path for general convex loss and penalty functions. We use the BLasso algorithm in our experimental section due to its ease of implementation and flexibility: the same code was used for different settings of the CAP penalty.

Similarly to boosting [10] and the Forward Stagewise Fitting algorithm [8], the BLasso algorithm works by taking forward steps of fixed size in the direction of steepest descent of the loss function. However, BLasso also allows for backward steps that take the penalty into account. With the addition of these backward steps, the BLasso is proven to approximate the Lasso path arbitrarily close, provided the step size can get small. An added advantage of the algorithm is its ability to trade off between precision and computational expense by adjusting the step size. For a detailed description of the algorithm, we refer the reader to Zhao and Yu [24].

3.1.2. *Piecewise linear paths: $L_2$-loss, $L_1$-norm and $L_\infty$-norm penalization.* For piecewise quadratic, convex loss functions and $\gamma_k \in \{1, \infty\}$ for all $k = 0, \ldots, K$, the CAP regularization path is known to be piecewise linear [18]. In these cases, it is possible to devise algorithms that jump from one breakpoint to the next while exactly computing their respective estimates as in the homotopy/LARS-LASSO algorithm [8, 17]. Next, we introduce two such algorithms for the $L_2$-loss: the first (iCAP) for grouped selection and the second for hierarchical selection (hiCAP). Before that, we present an algorithm for the $L_2$-loss estimates penalized by the $L_\infty$-norm (iLASSO). It serves as a stepping stone between the homotopy/LARS-LASSO [8, 17] and the iCAP and hiCAP path-tracing algorithms.

*The regularization path for iLASSO ($\infty$-LASSO).* The iLASSO estimate corresponds to the bridge regression (2) with the $L_2$-loss and $\gamma = \infty$. The KKT conditions defining the estimate for a particular $\lambda$ are

$$
\begin{bmatrix} \mathcal{X}'_{\mathcal{R}_\lambda} \mathcal{X}_{\mathcal{R}_\lambda} & \mathcal{X}_{\mathcal{R}_\lambda} X_{\mathcal{U}_\lambda} \\ X'_{\mathcal{U}_\lambda} \mathcal{X}_{\mathcal{R}_\lambda} & X'_{\mathcal{U}_\lambda} X_{\mathcal{U}_\lambda} \end{bmatrix} \begin{bmatrix} \hat{\alpha} \\ \hat{\beta}_{\mathcal{U}_\lambda} \end{bmatrix} = \begin{bmatrix} \mathcal{X}'_{\mathcal{R}_\lambda} Y - \lambda \\ X_{\mathcal{U}_\lambda} Y \end{bmatrix},
$$
(10)
$$\hat{\beta}_{\mathcal{R}_\lambda} = \hat{\alpha} S_\lambda,$$

where $\mathcal{R}_\lambda = \{j : |\hat{\beta}_j| = \|\hat{\beta}\|_\infty\}$, $\mathcal{U}_\lambda = \{j : |\hat{\beta}_j| < \|\hat{\beta}\|_\infty\}$, $S(\lambda) = \text{signs}[X'(Y - X\hat{\beta}(\lambda))]$, and $\mathcal{X}_{\mathcal{R}_\lambda} = \sum_{j \in \mathcal{R}_\lambda} S_j(\lambda) X_j$. From these conditions, it follows that $|\hat{\beta}_j(\lambda)| = \hat{\alpha}$, for all $j \in \mathcal{R}_\lambda$ and $X'_j(Y - X\hat{\beta}(\lambda)) = 0$, for all $j \in \mathcal{U}_\lambda$. Starting from a breakpoint $\lambda_0$ and its respective estimate $\hat{\beta}(\lambda_0)$, the path moves in a direction $\Delta \hat{\beta}$ that preserves the KKT conditions. The next breakpoint is then determined by $\hat{\beta}_{\lambda_1} = \hat{\beta}_{\lambda_0} + \delta \cdot \Delta \hat{\beta}$ where $\delta > 0$ is the least value to cause an index to move between the $\mathcal{R}_\lambda$ and $\mathcal{U}_\lambda$ sets. The pseudo-code for the iLASSO algorithm is presented in the technical report version of this paper [25]. We now extend this algorithm to handle the grouped case.

*The iCAP algorithm ($\infty$-CAP).* The iCAP algorithm is valid for the $L_2$-loss and nonoverlapping groups with $\gamma_0 = 1$ and $\gamma_k \equiv \infty$. This algorithm operates on two levels: it behaves as the Lasso on the group level and as the iLASSO within each group. To make this precise, first define the *$k$th group correlation* at $\lambda$ to be $c_k(\lambda) = \|X'_{\mathcal{G}_k}(Y - X\beta(\lambda))\|_1$ and the *set of active groups* $\mathcal{A}_\lambda = \{j \in \{1, \ldots, K\} : |c_j(\lambda)| = \max_{k=1,\ldots,K} |c_k(\lambda)|\}$. At a given $\lambda$, the groups not in $\mathcal{A}_\lambda$ have all their coefficients set to zero and $\hat{\beta}(\lambda)$ is such that all groups in $\mathcal{A}_\lambda$ have the same group correlation size. At the within-group level, for $\hat{\beta}(\lambda)$ to be a solution, each index $j \in \mathcal{G}_k$ must belong to either of two sets: $\mathcal{U}_{\lambda,k} = \{j \in \mathcal{G}_k : X'_j(Y - X\hat{\beta}(\lambda)) = 0\}$ or $\mathcal{R}_{\lambda,k} = \{j \in \mathcal{G}_k : \hat{\beta}_j(\lambda) = \|\hat{\beta}_{\mathcal{G}_k}(\lambda)\|_\infty\}$.

For $\lambda_0 = \max_{j \in \infty,\ldots,\mathcal{K}} \|X'_{\mathcal{G}_j} Y\|_1$, the solution is given by $\hat{\beta}(\lambda_0) = 0$. From this point, a direction $\Delta \hat{\beta}$ can be found so that the conditions above are satisfied by $\hat{\beta}(\lambda_0) + \delta\Delta\hat{\beta}$ for small enough $\delta > 0$. To find the next breakpoint, compute the least value of $\delta > 0$ that causes one of the following events: a group is added to or removed from $\mathcal{A}_\lambda$; an index moves between the sets $\mathcal{U}_{\lambda,k}$ and $\mathcal{R}_{\lambda,k}$ for some group $k$; or a sign change occurs in the correlation between the residuals and a variable in an inactive group. If no such $\delta > 0$ exists, the algorithm moves toward an un-regularized solution along the direction $\Delta\hat{\beta}$. The pseudo-code is given in Appendix B. The Matlab code implementing this algorithm can be downloaded from http://www.stat.berkeley.edu/twiki/Research/YuGroup/Software.

*The hiCAP algorithm (hierarchical-$\infty$-CAP).* We now introduce an algorithm for hierarchical selection. It is valid for the $L_2$-loss when $\gamma_0 = 1$, $\gamma_k \equiv \infty$ and the graph representing the hierarchy is a tree. The KKT conditions in this case are to an extent similar to those of the nonoverlapping groups case. The difference is that the groups now change dynamically along the path. Specifics of the algorithm are lengthy to describe and do not provide much statistical insight. We give a high level description here and refer readers interested in implementing the algorithm to the code available at http://www.stat.berkeley.edu/twiki/Research/YuGroup/Software.

The algorithm starts by forming nonoverlapping groups such that:

- Each nonoverlapping group consists of a sub-tree.
- Viewing each of the subtrees formed as a supernode, the derived supertree formed by these supernodes must satisfy the condition that average correlation size (unsigned) between $Y$ and $X$'s within the supernode is higher than that of all its descendant supernodes.

Once these groups are formed, the algorithm starts by moving the coefficients in the root group as in the nonoverlapping iCAP algorithm. The optimality conditions are met because the root group has the highest average correlation. Then, the algorithm proceeds observing two constraints:

1. Average unsigned correlation between $Y - X\hat{\beta}(\lambda)$ and $X$'s within each supernode is at least as high than that of all its descendant supernodes.
2. Maximum unsigned coefficient of each supernode is larger than or equal to that of any of its descendants.

Between breakpoints, the path is found by determining a direction such that these conditions are met. Breakpoints are found by noticing they are characterized by:

- If the average correlation between $Y - X\hat{\beta}(\lambda)$ and a subtree $\mathcal{G}_a$ contained by a supernode equals that of a supernode, then $\mathcal{G}_a$ splits into a new supernode.

- If a supernode $a$'s maximum coefficient size equals that of a descendant supernode $b$, then they are combined into a new supernode. These would also guarantee that a super node with all zero coefficients should have descendants with all zero coefficients.
- If a supernode with all zero coefficients and a descendant reached equal average correlation size (unsigned), they are merged.

3.2. *Choosing the regularization parameter $\lambda$.* For the selection of the regularization parameter $\lambda$ under general CAP penalties, we propose the use of cross-validation. We refer the reader to Stone [20] and Efron [6] for details on cross-validation. For the particular case of the iCAP having nonoverlapping groups and under the $L_2$-loss, we extend the Zou, Hastie and Tibshirani [27] unbiased estimate of the degrees of freedom for the LASSO. This estimate can then be used in conjunction with an information criterion to pick an appropriate level of regularization. We choose to use Sugiura's [21] $AIC_C$ information criterion. Our preference for this criterion follows from its being a correction of Akaike's [1] $AIC$ for small samples. Letting $k$ denote the dimensionality of a linear model, the $AIC_C$ criterion for a linear model is

$$(11) \qquad AIC_C = \frac{n}{2}\log\left(\sum_{i=1}^{n}(Y_i - X_i\hat{\beta}(\lambda))^2\right) + \frac{n}{2}\left(\frac{1+k/n}{1-k+2/n}\right),$$

where $n$ corresponds to the sample size. A model is selected by minimizing the expression in (11). To apply this criterion to iCAP, we substitute $k$ by estimates of the effective number of parameters (degrees of freedom). As we will see in our experimental section, the $AIC_C$ criterion resulted in good predictive performance, even in the "small-n, large-p" setting.

The extension of Zou, Hastie and Tibshirani [27] unbiased estimates of the degrees of freedom to iCAP follows from noticing that, between breakpoints, the LASSO, iLASSO and iCAP estimates mimic the behavior of projections on linear subspaces spanned by subsets of the predictor variables. It follows that standard results for linear estimates establish the dimension of the projecting subspace as an unbiased estimate of the degrees of freedom of penalized estimates in a broader set of penalties. Therefore, to compute unbiased estimates of degrees of freedom for the iLASSO or iCAP regressions above, it is enough to count the number of free parameters being fit at a point in the path. Letting $\mathcal{A}_\lambda$ and $\mathcal{U}_{k,\lambda}$ be as defined in Section 3.1 above, the resulting estimates of the degrees of freedom are $\hat{df}(\lambda) = |\mathcal{U}_\lambda| + 1$ for the iLASSO and $\hat{df}(\lambda) = |\mathcal{A}_\lambda| + \sum_{k\in\mathcal{A}_\lambda}|\mathcal{U}_{k,\lambda}|$ for the iCAP. A complete proof for the iLASSO case and a sketch of the proof for the iCAP case are given in Appendix A.1. Based on the above reasoning, we believe similar results should hold under the $L_2$-loss for any CAP penalties built exclusively from the $L_1$-norm and the $L_\infty$-norm. The missing ingredient is a way of determining the dimension of a projecting subspace along the path in broader settings. Future research will be devoted to that.

**4. Experimental results.** We now illustrate and evaluate the use of CAP in a series of simulated examples. The CAP framework needs the input of a grouping or hierarchical structure. For the grouping simulations, we explore the possibility of having a data-driven group structure choice and show that CAP enjoys a certain degree of robustness to misspecified groups. We leave a completely data-driven group structure choice as a topic of future research. For the hierarchical selection, two examples illustrate the use of the framework presented in Section 2.2.2 and the penalty in (9) to induce a pre-determined hierarchy. We defer the study of data-defined hierarchies for future work.

We will be comparing the predictive performance, sparsity and parsimony of different CAP estimates and that of the LASSO. As a measure of prediction performance we use the model error $ME(\hat{\beta}) = (\hat{\beta} - \beta)\mathbb{E}(X'X)(\hat{\beta} - \beta)$. In what concerns sparsity, we look not only at the number of variables selected by CAP and LASSO, but also at measures of sparsity that take the added structure into account: the number of selected groups for group selection and a properly defined measure for hierarchical selection (see Section 4.2). Whenever possible, we also compare the parsimony of the selected models as measured by the effective degrees of freedom and compare the cross-validation and $AIC_C$ based selections for the tuning parameter $\lambda$.

In all experiments below, the data is simulated from

$$Y = X\beta + \sigma\varepsilon, \qquad (12)$$

with $\varepsilon \sim \mathcal{N}(0, \mathbb{I})$. The parameters $\beta$, $\sigma$ as well as the covariance structure of $X$ are set specifically for each experiment.

4.1. *Grouped selection results.* In our grouping experiments, the group structure among the $p$ predictors in $X$ is due to their relationship to a set of $K$ zero-mean Gaussian hidden factors $Z \in \mathbb{R}^K$:

$$\mathrm{cov}(Z_k, Z_{k'}) = \begin{cases} 2.0, & \text{if } |k - k'| = 0, \\ 1.0, & \text{if } |k - k'| = 1, \\ 0.0, & \text{if } |k - k'| > 1, \end{cases} \qquad \text{for } k, k' \in \{1, \ldots, K\}.$$

A predictor $X_j$ in group $\mathcal{G}_k$ is the sum of the $Z_k$ factor and a noise term

$$X_j = Z_k + \eta_j \qquad \text{if } j \in \mathcal{G}_k,$$

with the Gaussian $\eta_j$ noise term having zero mean and

$$\mathrm{cov}(\eta_j, \eta_{j'}) = (4.0) \cdot 0.95^{|j-j'|}.$$

The correlation between the factors $Z$ is meant to make group selection slightly more challenging. Also, we chose the covariance structure of the disturbance terms $\eta$ to avoid an overly easy case for the empirical determination of the groups described below.

*Clustering for forming groups.* In applications, the grouping structure must often be estimated from data. In our grouping experiments, we estimate the true group structure $\mathcal{G}$ by clustering the predictors in $X$ using the Partitioning Around Medoids (PAM) algorithm (see [13]). The PAM algorithm needs the number of groups as an input. Instead of trying to fine-tune the selection of the number of groups, we fix the number of groups $\tilde{K}$ in the estimated clustering above, below and at the correct value $K$, specifically, we set:

- $\tilde{K} = K$: proper estimation of the number of groups;
- $\tilde{K} = 0.5 \cdot K$: severe underestimation of the number of groups;
- $\tilde{K} = 1.5 \cdot K$: severe overestimation of the number of groups.

Implicitly, we are assuming that a tuning method that can set the number of estimated groups $\tilde{K}$ in PAM to be between $0.5K$ (alt. below $1.5K$) and the true number of groups $K$ has results that are no worse than the ones observed in our underestimated (alt. overestimated) scenario below.

4.1.1. *Effect of group norms and $\lambda$ selection methods.* In this first experiment, we want to compare the difference among CAP estimates using alternative settings for the within-group norm and the LASSO. We keep the dimensionality of the problem low and emulate a high-dimensional setting ($n < p$) by setting

$$n = 80, \qquad p = 100 \quad \text{and} \quad K = 10.$$

The coefficients $\beta$ are made dissimilar (see Figure 6) within a group to avoid undue advantage to iCAP:

$$\beta_j = \begin{cases} 0.10(1 + 0.9^{j-1}), & \text{for } j = 1, \ldots, 10; \\ 0.04(1 + 0.9^{j-11}), & \text{for } j = 11, \ldots, 20; \\ 0.01(1 + 0.9^{j-21}), & \text{for } j = 21, \ldots, 30; \\ 0, & \text{otherwise.} \end{cases}$$

The noise level is set to $\sigma = 3$ and results are based on 50 replications.

The results reported in Table 1 show that all the different CAP penalties considered significantly reduced the model error in the comparison with the LASSO. The reduction in model error was also observed to be robust to misspecifications of the group structure (the cases $\tilde{K} = 0.5 \cdot K$ and $\tilde{K} = 1.5 \cdot K$).

Table 2 shows that our estimate for the degrees of freedom used in conjunction with the $AIC_C$ criterion was able to select predictive models as good as 10-fold cross-validation at a lower computational expense. The comparison of the results across the number of clusters used to group the predictors show that the improvement in prediction was robust to misspecifications in the number of groups used to cluster the predictors. iCAP's performance was the most sensitive to this type of misspecification as the $L_\infty$-norm makes a heavier use of the prespecified grouping information.

TABLE 1
**Comparison of different CAP penalties.** *Results for the LASSO and different CAP penalties for the first grouping experiment (based on 50 replications). All CAP penalties considered reduced the mean model error and number of selected groups in the comparison with the LASSO. The CAP penalty with $\gamma_k \equiv 4$ had a slight advantage over iCAP and the GLASSO*

| N. groups | LASSO | GLASSO | CAP(4) | iCAP |
|---|---|---|---|---|
| | | Model errors | | |
| Underestimated | 1.863 | 1.025 | **0.918** | 1.429 |
| (0.5K) | (0.194) | (0.101) | **(0.106)** | (0.316) |
| Right | 1.863 | 1.048 | **0.835** | 0.933 |
| (1.0K) | (0.194) | (0.094) | **(0.100)** | (0.092) |
| Overestimated | 1.863 | 1.159 | **0.970** | 1.271 |
| (1.5K) | (0.194) | (0.089) | **(0.090)** | (0.135) |
| | | Number of selected variables | | |
| Underestimated | **13.567** | 45.233 | 39.650 | 65.333 |
| (0.5K) | **(1.243)** | (3.504) | (3.490) | (4.633) |
| Right | **13.567** | 38.200 | 32.450 | 49.000 |
| (1.0K) | **(1.243)** | (2.502) | (1.748) | (3.567) |
| Overestimated | **13.567** | 33.900 | 33.450 | 48.400 |
| (1.5K) | **(1.243)** | (2.452) | (3.097) | (3.461) |
| | | Number of selected groups | | |
| Underestimated | 6.233 | 5.600 | **4.600** | 6.933 |
| (0.5K) | (0.491) | (0.400) | **(0.387)** | (0.442) |
| Right | 6.233 | 4.067 | **3.250** | 4.900 |
| (1.0K) | (0.491) | (0.275) | **(0.176)** | (0.357) |
| Overestimated | 6.233 | 4.100 | **3.950** | 5.333 |
| (1.5K) | (0.491) | (0.326) | **(0.352)** | (0.402) |

TABLE 2
**Comparison of model errors according to $\lambda$ selection method.** *Mean difference in model error between CAP and LASSO models selected using the $AIC_C$ criterion and 10-fold cross-validation. The $AIC_C$ criterion had results comparable to 10-fold cross-validation*

| | ME($AIC_C$)-ME(CV) | |
|---|---|---|
| | **LASSO ($\gamma_k \equiv 1$)** | **iCAP ($\gamma_k \equiv \infty$)** |
| Underestimated | $-0.253$ | $-0.077$ |
| (0.5K) | (0.177) | (0.048) |
| Right | $-0.470$ | $-0.267$ |
| (1.0K) | (0.388) | (0.207) |
| Overestimated | $-0.324$ | $-0.112$ |
| (1.5K) | (0.245) | (0.065) |

In terms of sparsity, the CAP estimates include a larger number of variables than the LASSO due to its block-inclusion nature. If we look at how many of the true groups are selected instead, we see that the CAP estimates made use of a lesser number of groups than the LASSO: an advantage if group selection is the goal. The low ratio between the number of variables and number of groups selected by the LASSO provide evidence that the LASSO estimates did not preserve the group structure selecting only a few variables from each group.

4.1.2. *Grouping with small-n-large-p.* We now compare iCAP and the LASSO when the number of predictors $p = Kq$ grows due to an increase in either the number of groups $K$ or in the group size $q$. The sample size is fixed at $n = 80$. The coefficients are randomly selected for each replication according to two different schemes: in the Grouped Laplacian scheme the coefficients are constant within each group and equal to $K$ independent samples from a Laplacian distribution with parameter $\alpha_G$; in the Individual Laplacian scheme the $p$ coefficients are independently sampled from a Laplacian distribution with parameter $\alpha_I$. The Grouped Laplacian scheme favors iCAP due to the grouped structure of the coefficients whereas the Individual Laplacian scheme favors the LASSO. The parameters $\alpha_G$ and $\alpha_I$ were adjusted so the signal power $\mathbb{E}(\beta'X'X\beta)$ is roughly constant across experiments. The complete set of parameters used is shown in Table 3.

We only report the results obtained from using the $AIC_C$ criterion to select $\lambda$. The results for 10-fold cross-validation were similar. Tables 4 and 5 show the results based on 100 replications.

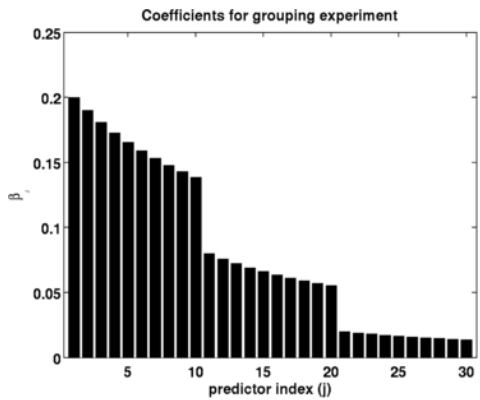

FIG. 6. **Profile of coefficients for grouping experiment.** *In the first grouping experiment, only the first three groups have nonzero coefficients. Within each group, the coefficients have an exponential decay.*

GROUPED AND HIERARCHICAL SELECTION 19The iCAP estimates had predictive performance better than or comparable to LASSO estimates. For the Grouped Laplacian case, the reduction in model error was very pronounced for all settings considered. In the Individual Laplacian case, iCAP and LASSO had comparable model errors for $p = 100$. As the number of predictors increased, however, iCAP resulted in lower model errors than the LASSO even under the Individual Laplacian regime. This result provides evidence that collecting highly correlated predictors into groups is beneficial in terms of predicting performance, especially as the ratio $p/n$ increases. The predictive gains from using CAP were also preserved under misspecified groupings. iCAP estimates had smaller or comparable model errors than the LASSO when the predictors were clustered into a number of groups that was 50% smaller or larger than the actual number of predictor clusters.

Contrary to what was expected, the iCAP estimates involved a number of groups comparable to the LASSO in all simulated scenarios. We believe this is due to a more parsimonious use of degrees of freedom by iCAP estimates. As more groups are added to the model, the within-group restrictions imposed by the $L_\infty$-norm prevent a rapid increase in the use of degrees of freedom. As a result, iCAP estimates can afford to include groups in the model more aggressively in an attempt to reduce the $L_2$-loss. This reasoning is supported by the lesser degrees of freedom used by iCAP selected models in the comparison to the LASSO.

4.2. *Hierarchical selection results.* Now, we provide examples where the hierarchical structure of the predictors is exploited to enhance the predictive performance of models. We define *hierarchical gap* as a measure of compliance to the given hierarchical structure: it is the minimal number of additional variables that should be added to the model so the hierarchy is satisfied. If a model satisfies a given hierarchy, no additional variables must be added and this measure equals zero.

4.2.1. *Hierarchical selection for ANOVA model with interaction terms.* We now further develop the regression model with interaction terms introduced in Section 2 (cf. Figure 4).

Table 3
*Parameters for the simulation of the small-n-large-p case*

| | | | Grouped Laplacian | | | Individual Laplacian | | | |
|---|---|---|---|---|---|---|---|---|---|
| $p$ | $q$ | $K$ | $\alpha_G$ | $\mathbb{E}(\beta'X'X\beta)$ | SNR | $\alpha_I$ | $\mathbb{E}(\beta'X'X\beta)$ | SNR | $\sigma$ |
| 100 | 10 | 10 | $1.00 \cdot 10^{-1}$ | 54.02 | 3.95 | $3.00 \cdot 10^{-1}$ | 54.00 | 3.95 | 3.7 |
| 250 | 10 | 25 | $0.63 \cdot 10^{-1}$ | 53.60 | 3.92 | $1.90 \cdot 10^{-1}$ | 54.15 | 3.96 | 3.7 |
| 250 | 25 | 10 | $0.43 \cdot 10^{-1}$ | 54.61 | 3.99 | $1.90 \cdot 10^{-1}$ | 54.15 | 3.96 | 3.7 |

TABLE 4
**Results for small-n-large-p experiment under Grouped Laplacian sampling.** *Results based on* 100 *replications and $AIC_C$ selected $\lambda$. The true model has its parameters sampled according to the Grouped Laplacian scheme (see Section [4.1.2](#)). The inclusion of grouping structure improves the predictive performance of the models whether the predictors are clustered in the correct number of groups (1.0K) or not (0.5K and 1.5K). LASSO selects a smaller number of variables and is slightly sparser in terms of number of groups. iCAP estimates are more parsimonious in terms of degrees of freedom*

|  | **LASSO** | **0.5K** | **1.0K** | **1.5K** |
|---|---|---|---|---|
|  | Model errors ||||
| $p=100$, $q=10$ | 5.028 | 3.783 | **2.839** | 3.481 |
|  | (0.208) | (0.172) | **(0.119)** | (0.132) |
| $p=250$, $q=10$ | 13.061 | 11.135 | **6.660** | 8.128 |
|  | (0.506) | (0.834) | **(0.227)** | (0.271) |
| $p=250$, $q=25$ | 8.356 | 5.113 | **3.479** | 4.457 |
|  | (0.379) | (0.228) | **(0.149)** | (0.202) |
|  | Number of selected variables ||||
| $p=100$, $q=10$ | **19.590** | 93.940 | 82.200 | 73.720 |
|  | **(0.546)** | (1.108) | (1.508) | (1.438) |
| $p=250$, $q=10$ | **26.070** | 211.090 | 169.700 | 144.740 |
|  | **(0.668)** | (3.335) | (2.866) | (2.781) |
| $p=250$, $q=25$ | **25.450** | 236.300 | 218.500 | 192.730 |
|  | **(0.664)** | (2.979) | (3.049) | (3.667) |
|  | Number of true groups selected ||||
| $p=100$, $q=10$ | **8.140** | 9.520 | 8.220 | 8.240 |
|  | **(0.163)** | (0.090) | (0.151) | (0.148) |
| $p=250$, $q=10$ | **14.550** | 21.780 | 16.970 | 16.720 |
|  | **(0.336)** | (0.290) | (0.287) | (0.317) |
| $p=250$, $q=25$ | **7.980** | 9.560 | 8.740 | 8.530 |
|  | **(0.160)** | (0.102) | (0.122) | (0.149) |
|  | Degrees of freedom ||||
| $p=100$, $q=10$ | 19.590 | 13.600 | **10.460** | 12.060 |
|  | (0.546) | (0.459) | **(0.365)** | (0.404) |
| $p=250$, $q=10$ | 26.070 | 19.710 | **18.490** | 20.190 |
|  | (0.668) | (0.646) | **(0.431)** | (0.501) |
| $p=250$, $q=25$ | 25.450 | 18.010 | **13.590** | 14.920 |
|  | (0.664) | (0.604) | **(0.481)** | (0.491) |

The data is generated according to (12) and set the number of variables to $d = 10$ resulting in a regression involving 10 main effects and 45 interactions:

$$X = [Z_1, \ldots, Z_{10}, Z_1 Z_2, \ldots, Z_1 Z_{10}, \ldots, Z_9 Z_{10}],$$
$$Z_j \stackrel{\text{i.i.d.}}{\sim} \mathcal{N}(0, 1).$$

GROUPED AND HIERARCHICAL SELECTION 21TABLE 5
**Results for small-n-large-p experiment under Individual Laplacian sampling.** *Results based on* 100 *replications and* $AIC_C$ *selected* $\lambda$. *The true model has its parameters sampled according to the Independent Laplacian scheme (see Section 4.1.2). For a lower-dimensional model* ($p = 100, q = 10$), *the predictive performance of iCAP is comparable to the LASSO. For higher dimensions* ($p = 250$, $q = 10$ *and* $q = 25$), *iCAP has better predictive performance. LASSO selects a smaller number of variables than iCAP and a comparable number of groups in all cases. iCAP estimates are still more parsimonious in terms of degrees of freedom*

|  | LASSO | 0.5K | 1.0K | 1.5K |
|---|---|---|---|---|
|  | Model errors | | | |
| $p = 100, q = 10$ | 10.310 | 10.885 | **10.153** | 11.056 |
|  | (0.309) | (0.388) | **(0.300)** | (0.348) |
| $p = 250, q = 10$ | 22.560 | 18.790 | **18.194** | 18.990 |
|  | (0.701) | (0.446) | **(0.463)** | (0.422) |
| $p = 250, q = 25$ | 19.891 | 17.483 | **16.544** | 18.301 |
|  | (0.614) | (0.424) | **(0.387)** | (0.493) |
|  | Number of selected variables | | | |
| $p = 100, q = 10$ | **20.200** | 96.460 | 90.700 | 73.530 |
|  | **(0.620)** | (0.869) | (1.249) | (1.459) |
| $p = 250, q = 10$ | **25.540** | 228.070 | 180.700 | 150.000 |
|  | **(0.620)** | (2.246) | (3.003) | (2.508) |
| $p = 250, q = 25$ | **24.440** | 243.490 | 234.250 | 198.040 |
|  | **(0.589)** | (1.530) | (1.935) | (2.638) |
|  | Number of true groups selected | | | |
| $p = 100, q = 10$ | 9.580 | 9.760 | 9.480 | **9.320** |
|  | (0.106) | (0.092) | (0.098) | **(0.119)** |
| $p = 250, q = 10$ | 21.210 | 24.110 | 21.580 | **20.720** |
|  | (0.309) | (0.115) | (0.266) | **(0.290)** |
| $p = 250, q = 25$ | 9.720 | 9.870 | 9.730 | **9.650** |
|  | (0.057) | (0.049) | (0.066) | **(0.069)** |
|  | Degrees of freedom | | | |
| $p = 100, q = 10$ | 20.200 | 16.320 | 16.140 | **14.650** |
|  | (0.620) | (0.684) | (0.635) | **(0.614)** |
| $p = 250, q = 10$ | 25.540 | 22.230 | **20.290** | 21.210 |
|  | (0.620) | (0.630) | **(0.486)** | (0.536) |
| $p = 250, q = 25$ | 24.440 | 20.360 | 20.060 | **18.670** |
|  | (0.589) | (0.610) | (0.635) | **(0.637)** |

We assume the hierarchical structure is given that the second order terms are to be added to the model only after their corresponding main effects. This is an extension of the hierarchical structure in Figure 4 from $d = 4$ to

TABLE 6
*Simulation setup for the ANOVA experiment*

|  | **Coefficients** | | | | | | | |
| --- | --- | --- | --- | --- | --- | --- | --- | --- |
| Description | $Z_1Z_2$ | $Z_1Z_3$ | $Z_1Z_4$ | $Z_2Z_3$ | $Z_2Z_4$ | $Z_3Z_4$ | $\sigma$ | SNR |
| No interactions | 0 | 0 | 0 | 0 | 0 | 0 | 3.7 | 4.02 |
| Weak interactions | 0.5 | 0 | 0 | 0.1 | 0.1 | 0 | 3.7 | 4.08 |
| Moderate interactions | 1.0 | 0 | 0 | 0.5 | 0.4 | 0.1 | 3.7 | 4.33 |
| Strong interactions | 5 | 0 | 0 | 4 | 2 | 0 | 3.7 | 13.88 |
| Very strong interactions | 7 | 7 | 7 | 2 | 2 | 1 | 3.7 | 38.20 |

$d = 10$. Applying (9) to this graph with uniform $\alpha_m$ weights gives:

$$T(\beta) = \sum_{j=1}^{10} \sum_{i=1}^{j-1} [|\beta_{i,j}| + \|(\beta_i, \beta_j, \beta_{i,j})\|_{\gamma_{i,j}}].$$

In this case, each interaction term is penalized in three factors of the summation, which agrees to the number of variables that are added to the model ($Z_{ij}$, $Z_i$ and $Z_j$).

We set the first four main effect coefficients to be $\beta_1 = 7$, $\beta_2 = 2$, $\beta_3 = 1$ and $\beta_4 = 1$ with all remaining main effects set to zero. The main effects are kept fixed throughout and five different levels of interaction strengths are considered as shown in Table 6. The variance of the noise term was kept constant across the experiments. The number of observations $n$ was set to 121. The results reported in Table 7 were obtained from 100 replications.

For low and moderate interactions, the introduction of hierarchical structure reduced the model error from the LASSO. For strong interactions, the CAP and LASSO results were comparable. For very strong interactions, the implicit assumption of smaller second order effects embedded in the hierarchical selection is no longer suitable causing the LASSO to reach a better predictive performance.

In addition, CAP selected models involving on average a slightly lesser or equal number of variables to the LASSO in all simulated cases. The hierarchical gap (see definition in Section 4.2) for the LASSO shows that it did not comply with the hierarchy of the problem. According to the theory developed in Section 2, for CAP estimates this difference should be exactly zero. The small deviations from zero observed in our experiments are due to the approximate nature of the BLasso algorithm.

4.2.2. *Multiresolution model.* In this experiment, the true signal is given by a linear combination of Haar wavelets at different resolution levels. Let-

TABLE 7
**Simulation results for the hierarchical ANOVA example.** *Results based on* 50 *replications,* 121 *observations and 10-fold CV. Hierarchical structure lead to reduced model error and sparser models. The hierarchy gap is the number of variables that must be added to the model so the hierarchy is satisfied. The LASSO does not respect the model hierarchy. The small deviations from zero for CAP estimates are due to BLASSO approximation*

|  | **LASSO** | **"GLASSO"** | **CAP(4)** | **iCAP** |
|---|---|---|---|---|
| | | Model errors | | |
| No interactions | 3.367 | 1.481 | 1.478 | **1.466** |
| | (0.288) | (0.133) | (0.134) | **(0.124)** |
| Weak interactions | 4.032 | 2.190 | 2.296 | **2.117** |
| | (0.303) | (0.147) | (0.161) | **(0.117)** |
| Moderate interactions | 5.905 | 4.260 | 4.090 | **4.085** |
| | (0.307) | (0.205) | (0.207) | **(0.196)** |
| Strong interactions | 8.912 | 8.901 | **7.793** | 8.388 |
| | (0.695) | (0.621) | **(0.568)** | (0.626) |
| Very strong int. | **11.474** | 11.998 | 14.538 | 26.072 |
| | **(0.758)** | (0.800) | (0.915) | (1.351) |
| | | Number of selected variables | | |
| No interactions | 13.440 | 12.720 | 11.440 | **9.280** |
| | (0.897) | (0.935) | (0.725) | **(0.431)** |
| Weak interactions | 13.960 | 13.780 | 12.720 | **9.920** |
| | (0.936) | (1.041) | (0.806) | **(0.442)** |
| Moderate interactions | 16.160 | 17.240 | 14.960 | **11.260** |
| | (1.021) | (1.015) | (0.773) | **(0.488)** |
| Strong interactions | 21.800 | 27.000 | 20.680 | **13.740** |
| | (0.965) | (0.910) | (0.636) | **(0.384)** |
| Very strong int. | 26.400 | 28.020 | 20.520 | **14.560** |
| | (0.670) | (0.577) | (0.448) | **(0.289)** |
| | | Hierarchy gap | | |
| No interactions | 3.560 | **0.220** | 0.420 | 0.800 |
| | (0.192) | **(0.066)** | (0.107) | (0.128) |
| Weak interactions | 3.560 | **0.160** | 0.340 | 0.840 |
| | (0.169) | **(0.052)** | (0.079) | (0.112) |
| Moderate interactions | 3.440 | **0.420** | 0.640 | 1.020 |
| | (0.174) | **(0.091)** | (0.106) | (0.150) |
| Strong interactions | 4.240 | **0.740** | 1.600 | 2.020 |
| | (0.184) | **(0.102)** | (0.143) | (0.163) |
| Very strong int. | 3.780 | **1.080** | 1.580 | 0.920 |
| | (0.155) | **(0.140)** | (0.143) | (0.137) |

ting $Z_{ij}$ denote the Haar wavelet at the $j$th position of level $i$, we have

$$Z_{ij}(t) = \begin{cases} -1, & \text{if } t \in \left(\dfrac{j}{2^{i+1}}, \dfrac{j+1}{2^{i+1}}\right), \\ 1, & \text{if } t \in \left(\dfrac{j+1}{2^{i+1}}, \dfrac{j+2}{2^{i+1}}\right), \\ 0, & \text{otherwise}, \end{cases} \quad \text{for } i \in \mathbb{N}, j = 0, 1, 2, \ldots, 2^{i-1}.$$

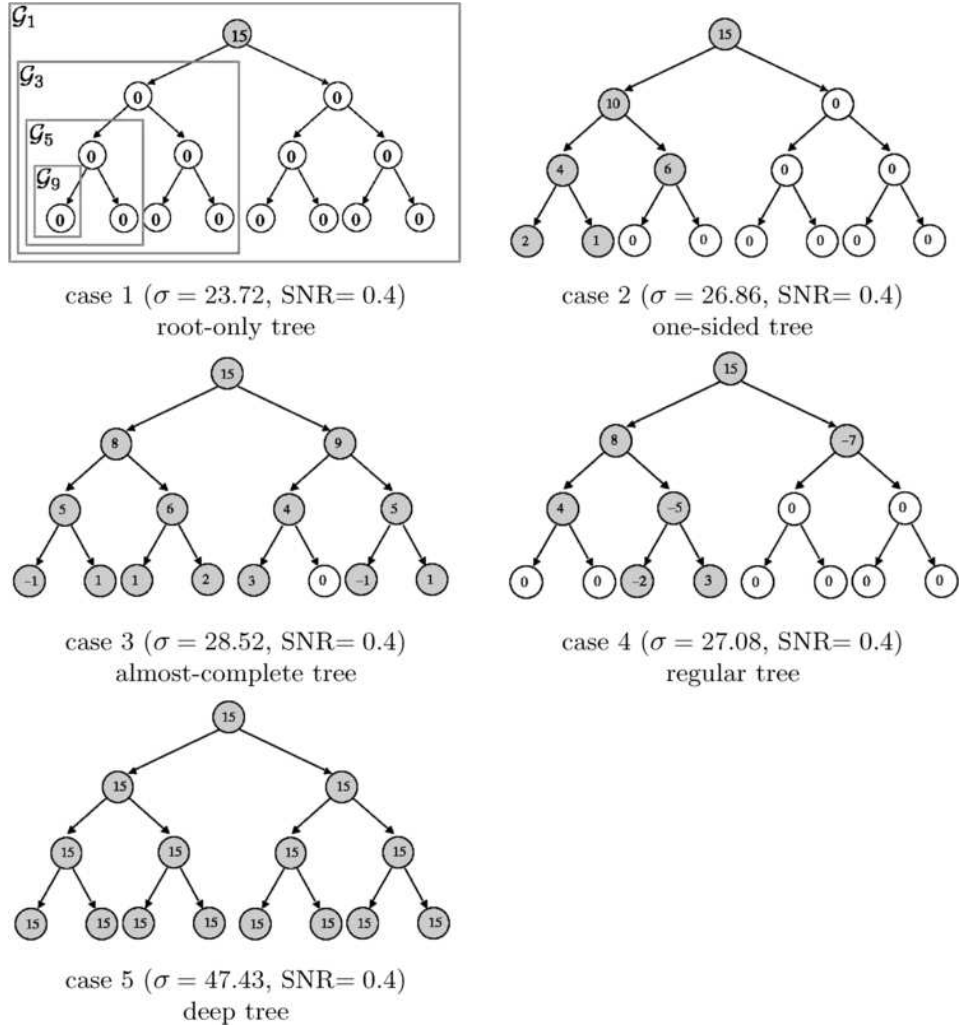

Fig. 7. *Coefficients and hierarchies used in the wavelet tree example.*

In what follows, we let a vector $\tilde{Z}_{ij}$ be formed by computing $Z_{ij}(t)$ at 16 equally spaced "time points" over $[0, 1]$. The design matrix $X$ is then given by $[\tilde{Z}_{00}, \quad \tilde{Z}_{10} \quad \tilde{Z}_{11}, \quad \tilde{Z}_{20} \quad \cdots \quad \tilde{Z}_{23}, \quad \tilde{Z}_{30} \quad \cdots \quad \tilde{Z}_{37}]$.

We consider five different settings for the value of $\beta$ with various sparsity levels and tree depths. The parameter $\sigma$ is adjusted to keep the signal to noise ratio at 0.4 in all cases. The true model parameters are shown within the tree hierarchy in Figure 7.

The simulated data corresponds to 5 sets of observations of the 16 "time points." Five-fold cross-validation was used for selecting the regularization level. At each cross-validation round, the 16 points kept out of the sample correspond to each of the "time positions" ("balanced cross-validation").

In its upper left panel, Figure 7 shows the directed graph representing the tree hierarchy used to form the CAP penalty using the recipe laid out by (9) in Section 2 with $\alpha_m = 1$ for all $m$. Our option for setting all weights to 1 in this case is similar to the one presented in the ANOVA experiment above: the number of variables added to the model when a variable is added matches the number of times its coefficient appears on the penalty.

The results for the hierarchical selection are shown in Table 8. The use of the hierarchical information greatly reduced the model error as well as the number of selected variables. As in the ANOVA cases, the hierarchical gap shows that the LASSO models do not satisfy the tree hierarchy. The approximate nature of the BLASSO algorithm again causes the hierarchical gap of GLASSO and CAP(4) to deviate slightly from zero. For hiCAP estimates, perfect agreement with the hierarchy is observed as the hiCAP algorithm (see Section 3) is exact.

4.3. *Additional experiments.* In addition to the experiments presented above, we have also run CAP under examples taken from Yuan and Lin [23] and Zou and Hastie [26]. The results are similar to the ones obtained above: CAP results in improved prediction performance over the LASSO with models involving a larger number of variables, a similar or smaller number of groups and making use of less degrees of freedom. We invite the reader to check the details in the technical report version of this paper [25].

**5. Discussion and concluding remarks.** In this paper, we have introduced the Composite Absolute Penalty (CAP) family. It provides a regularization framework for incorporating predetermined grouping and hierarchical structures among the predictors by combining $L_\gamma$-norm penalties and applying them to properly defined groups of coefficients.

The definition of the groups to which norms are applied is instrumental in determining the properties of CAP estimates. Nonoverlapping groups give rise to group selection as done previously by Yuan and Lin [23] and Kim, Kim and Kim [14] where $L_1$ and $L_2$ norms are combined. CAP penalties extend these works by letting norms other than the $L_2$-norm to be applied to the groups of variables. Combinations of the $L_1$ and $L_\infty$ are convenient from a computational standpoint as illustrated by the iCAP (nonoverlapping groups) with fast homotopy/LARS-type algorithms. Its Matlab code can be downloaded from our research group website http://www.stat.berkeley.edu/twiki/Research/YuGroup/Software.

The definition of CAP penalties also generalizes previous work by allowing the groups to overlap. Here, we have shown how to construct overlapping groups leading to hierarchical selection. Combinations of the $L_1$ and $L_\infty$ are also computationally convenient to hierarchical selection as illustrated by the hiCAP algorithm for tree hierarchical selection (also available from our research group website).

TABLE 8
**Simulation results for the hierarchical wavelet tree example.** *Results based on* 200 *replications,* $5 \times 16$ *observations and 5 fold "balanced" CV. Hierarchical structure lead to reduced model error and sparser models. The hierarchy gap is the number of variables that must be added to the model so the hierarchy is satisfied. The LASSO does not respect the tree hierarchy. Small discrepancies in GLASSO and CAP(4) due to approximation in BLASSO*

|  | **LASSO** | **GLASSO** | **CAP(4)** | **iCAP** |
|---|---|---|---|---|
|  | Model errors | | | |
| Root-only tree | 40.508 | **26.994** | 28.498 | 28.909 |
|  | (2.039) | **(1.635)** | (1.650) | (1.820) |
| One-sided | 80.197 | **57.100** | 58.228 | 57.195 |
|  | (2.954) | **(2.519)** | (2.465) | (2.532) |
| Complete tree | 112.578 | **76.911** | 79.498 | 78.117 |
|  | (3.407) | **(2.661)** | (2.736) | (2.733) |
| Regular tree | 82.979 | **58.020** | 60.037 | 60.259 |
|  | (2.738) | **(2.237)** | (2.252) | (2.379) |
| Heavy-leaved tree | 454.607 | 388.262 | 385.770 | **359.154** |
|  | (11.950) | (10.015) | (9.884) | **(10.565)** |
|  | Number of selected variables | | | |
| Root-only tree | 4.070 | 4.495 | 4.405 | **3.185** |
|  | (0.209) | (0.265) | (0.241) | (0.244) |
| One-sided | 6.080 | 6.770 | 6.235 | **5.415** |
|  | (0.211) | (0.247) | (0.218) | (0.245) |
| Complete tree | 7.010 | 7.605 | 6.995 | **6.720** |
|  | (0.226) | (0.227) | (0.213) | (0.248) |
| Regular tree | 6.140 | 6.690 | 6.255 | **5.490** |
|  | (0.230) | (0.243) | (0.218) | (0.250) |
| Heavy-leaved tree | 10.985 | 11.630 | **10.930** | 11.240 |
|  | (0.258) | (0.192) | (0.186) | (0.239) |
|  | Hierarchy gap | | | |
| Root-only tree | 1.765 | 0.235 | 0.360 | **0.000** |
|  | (0.096) | (0.038) | (0.049) | **(0.000)** |
| One-sided | 0.710 | 0.180 | 0.300 | **0.000** |
|  | (0.064) | (0.031) | (0.038) | **(0.000)** |
| Complete tree | 1.640 | 0.310 | 0.505 | **0.000** |
|  | (0.091) | (0.042) | (0.057) | **(0.000)** |
| Regular tree | 1.210 | 0.225 | 0.335 | **0.000** |
|  | (0.069) | (0.030) | (0.039) | **(0.000)** |
| Heavy-leaved tree | 1.455 | 0.210 | 0.370 | **0.000** |
|  | (0.079) | (0.034) | (0.044) | **(0.000)** |

In a set of simulated examples, we have shown that CAP estimates using a given grouping or hierarchical structure can reduce the model error when compared to LASSO estimates. In the grouped-selection case, we show such reduction has taken place in the "small-$n$-large-$p$" setting and was observed even when the groups were data-determined (noisy) and the resulting number of groups was within a large margin of the actual number of groups among the predictors (between 50% and 150% of the true number of groups).

In addition, iCAP predictions are more parsimonious in terms of use of degrees of freedom being less sensitive to disturbances in the observed data [7]. Finally, CAP's ability to select models respecting the group and hierarchical structure of the problems makes its estimates more interpretable. It is a topic of our future research to explore ways to estimate group and hierarchical structures completely based on data.

## APPENDIX A: PROOFS

PROOF OF THEOREM 1. Algebraically, we have for $\gamma > 1$

$$\frac{\partial}{\partial \beta_1} T(\beta) = \frac{\partial}{\partial \beta_1} \|\beta\|_\gamma = \text{sign}(\beta_1) \left( \frac{|\beta_1|}{\|\beta\|_\gamma} \right)^{(\gamma-1)}.$$

As a result, if $\beta_2 > 0$, $\beta_1$ is locally not penalized at 0 and it will only stay at this point if the gradient of the loss function $L$ is exactly zero for $\beta_1 = 0$. Unless the distribution of the gradient of the loss function has an atom at zero for $\beta_1$, $\beta_1 \neq 0$ with probability one. □

PROOF OF THEOREM 2. It is enough to prove that $T$ is convex:

1. $\mathbf{T}(\alpha \cdot \beta) = |\alpha| \cdot \mathbf{T}(\beta)$, for all $\alpha \in \mathbb{R}$: For each group $k$, $N_k(\alpha\beta) = \alpha N_k(\beta)$. Thus $T(\alpha\beta) = \|\mathbf{N}(\alpha\beta)\|_{\gamma_0} = |\alpha| \|\mathbf{N}(\beta)\|_{\gamma_0} = |\alpha| T(\beta)$;
2. $\mathbf{T}(\beta_1 + \beta_2) \leq \mathbf{T}(\beta_1) + \mathbf{T}(\beta_2)$: Using the triangular inequality,

$$T(\beta_1 + \beta_2) = \sum_k (N_k(\beta_1 + \beta_2))^{\gamma_0} \leq \sum_k (N_k(\beta_1) + N_k(\beta_2))^{\gamma_0}$$
$$= \|N(\beta_1) + N(\beta_2)\|_{\gamma_0} \leq \|N(\beta_1)\|_{\gamma_0} + \|N(\beta_2)\|_{\gamma_0}$$
$$= T(\beta_1) + T(\beta_2).$$

Convexity follows by setting $\beta_1 = \theta\beta_3$ and $\beta_2 = (1-\theta)\beta_4$ with $\theta \in [0, 1]$. □

**A.1. DF Estimates for iLASSO and iCAP.** We now derive an unbiased estimate for the degrees of freedom of iLASSO fits along the regularization path. The optimization problem defining the iLASSO estimate is dual to the LASSO problem (see [2]). Facts 1 through 3 below follow from this duality and the results in Efron et al. [8] and Zou, Hastie and Tibshirani [27]. For the remainder of this section, we denote the iLASSO and iCAP fit by $\hat{\mu}(\lambda, y) = X\hat{\beta}(\lambda, y)$.

FACT A.1. *For each $\lambda$, there exists a set $\mathcal{K}_\lambda$ such that $\mathcal{K}_\lambda$ is the union of a finite collection of hyperplanes and for all $Y \in \mathcal{C}_\lambda = \mathbb{R}^n - \mathcal{K}_\lambda$, $\lambda$ is not a breakpoint in the regularization path.*

FACT A.2. *$\hat{\beta}(\lambda, y)$ is a continuous function of $y$ for all $\lambda$.*

FACT A.3. *If $y \in \mathcal{C}_\lambda$, then the sets $\mathcal{R}_\lambda$ and $\mathcal{U}_\lambda$ are locally invariant.*

From these three facts, we can prove:

LEMMA A.1. *For a fixed $\lambda \geq 0$ and $Y \in \mathcal{C}_\lambda$, $\hat{\mu}(\lambda, y)$ satisfies*

$$\|\hat{\mu}(\lambda, y + \Delta y) - \hat{\mu}(\lambda, y)\| \leq \|\Delta y\|, \qquad \text{for sufficiently small } \Delta y$$

*and*

$$\nabla \cdot \hat{\mu}(\lambda, y) = |\mathcal{U}_\lambda| + 1.$$

PROOF. We first notice that, $\hat{\mu}(\lambda, y) = X\hat{\beta}(\lambda, Y) = [\mathcal{X}_{\mathcal{R}_\lambda} \quad X_{\mathcal{U}_\lambda}] \cdot [\hat{\alpha} \quad \hat{\beta}'_{\mathcal{U}_\lambda}]'$. From the optimality conditions for the $L_\infty$ penalty,

$$(\mathcal{X}'\mathcal{X})\hat{\alpha}(\lambda, Y) = \mathcal{X}'Y - \lambda \cdot \text{sign}(Y - \mathcal{X}\hat{\alpha}(\lambda, Y))$$

and

$$(\mathcal{X}'\mathcal{X})\hat{\alpha}(\lambda, Y + \Delta Y) = \mathcal{X}'(Y + \Delta Y) - \lambda \cdot \text{sign}(Y + \Delta Y - \mathcal{X}\hat{\alpha}(\lambda, Y + \Delta Y)).$$

For $Y \in \mathcal{C}_\lambda$, there exists small $\Delta Y$ so the signs of the correlation between the residuals and each predictor are preserved. Subtracting the two equations above: $\hat{\mu}(\lambda, Y + \Delta Y) - \hat{\mu}(\lambda, Y) = \mathcal{X}(\mathcal{X}'\mathcal{X})^{-1}\mathcal{X}'\Delta Y$. Thus, $\hat{\mu}(\lambda, Y)$ behaves locally as a projection on a fixed subspace given by $\mathcal{R}_\lambda$ and $\mathcal{U}_\lambda$. From standard projection matrix results: $\|X\hat{\beta}(\lambda, y + \Delta y) - X\hat{\beta}(\lambda, y)\| \leq \|\Delta y\|$, for small $\Delta y$ and $\nabla \cdot \hat{\mu}(\lambda, Y) = \text{tr}(\mathcal{X}(\mathcal{X}'\mathcal{X})^{-1}\mathcal{X}') = |\mathcal{U}_\lambda| + 1$. □

Lemma A.1 implies that the fit $\hat{\mu}(\lambda, y)$ is uniformly Lipschitz on $\mathbb{R}^n$ (it is the closure of $\mathcal{C}_\lambda$). Using Stein's lemma and the divergent expression above:

THEOREM A.1. *The $L_\infty$-penalized fit $\hat{\mu}_\lambda(y)$ is uniformly Lipschitz for all $\lambda$. The degrees of freedom of $\hat{\mu}_\lambda(y)$ is given by $df(\lambda) = E[|\mathcal{U}_\lambda|] + 1$.*

The proof for the case of nonoverlapping groups follows the same steps. We only present a sketch of the proof as a detailed proof is not very insightful. Fact A.1 is proven by noticing that, for fixed $\lambda$, each of the conditions defining breakpoints require $Y$ to belong to a finite union of hyperplanes. Fact A.2 follows from the CAP objective function being convex and continuous in both $\lambda$ and $Y$. Fact A.3 is established by noticing that the sets $\mathcal{A}_\lambda$ and $\mathcal{R}_{k,\lambda}, \forall k = 1, \ldots, K$ are invariant in between breakpoints. As before, the CAP fit behaves (except for a shrinkage factor) as a projection onto a subspace whose dimension is the number of "free" parameters at that point of the path. The result follows from standard arguments for linear estimates.

## APPENDIX B: PSEUDO-CODE FOR THE ICAP ALGORITHM

1. Set $t = 0$, $\lambda_t = \max_k \|c_k(0)\|$, $\hat{\beta}(\lambda_t) = 0$.
2. Repeat until $\lambda_t = 0$:

   (a) Set $\mathcal{A}_{\lambda_t}$ to contain all groups with $c_k = \lambda_t$;
   
   (b) For each group $k$, set:
   
   $$\mathcal{U}_{\lambda,k} = \{j \in \mathcal{G}_k : X'_j(Y - X\hat{\beta}(\lambda_t)) = 0\} \quad \text{and} \quad \mathcal{R}_{\lambda,k} = \mathcal{G}_k - \mathcal{U}_{k,\lambda};$$
   
   (c) Determine a direction $\Delta\hat{\beta}$ such that:
   
   (i) if $k \notin \mathcal{A}_\lambda$, then $\Delta\hat{\beta}_{\mathcal{G}_k} = 0$;
   
   (ii) for $k \in \mathcal{A}_\lambda$, $\Delta\hat{\beta}_{\mathcal{R}_{k,\lambda}} = \alpha_k \cdot S_{\lambda,k}$ with $\alpha_k$ chosen so:
   
   $$c_k(\hat{\beta}(\lambda) + \delta \cdot \Delta\hat{\beta})$$
   $$= c_{k^*}(\hat{\beta}(\lambda) + \delta \cdot \Delta\hat{\beta}) \quad \text{for all } k, k^* \in \mathcal{A}_\lambda \quad \text{and}$$
   $$X'_{\mathcal{U}_{k,\lambda}}(Y - X(\hat{\beta}(\lambda) + \delta \cdot \Delta\hat{\beta}))$$
   $$= 0 \quad \text{for small enough } \delta > 0.$$
   
   (d) Compute the step sizes for which breakpoints occur:

$$\delta_A = \inf_{\delta > 0}\{c_{k^*}(\hat{\beta}_\lambda + \delta \cdot \Delta\hat{\beta}) = c_k(\hat{\beta}_\lambda + \delta \cdot \Delta\hat{\beta}) \text{ for some } k^* \notin \mathcal{A}_\lambda \text{ and } k \in \mathcal{A}_\lambda\},$$

$$\delta_I = \inf_{\delta > 0}\{\|\hat{\beta}_{\mathcal{G}_k}(\lambda) + \delta \cdot \Delta\hat{\beta}\|_\infty = 0 \text{ for some } k \in \mathcal{A}_\lambda\},$$

$$\delta_U = \inf_{\delta > 0}\{X'_m(Y - X(\hat{\beta}(\lambda) + \delta \cdot \Delta\hat{\beta})) = 0 \text{ for some } m \in \mathcal{G}_k \text{ with } k \in \mathcal{A}_\lambda\},$$

$$\delta_R = \inf_{\delta > 0}\{|\hat{\beta}_m(\lambda) + \delta \cdot \Delta\hat{\beta}_m| = \|\hat{\beta}_{\mathcal{G}_k}(\lambda) + \delta \cdot \Delta\hat{\beta}\|_\infty$$

for some $m \in \mathcal{U}_{k,\lambda}$ with $k \in \mathcal{A}_\lambda\}$,

$$\delta_S = \inf_{\delta > 0} \{X'_m(Y - X(\hat{\beta}(\lambda) + \delta \cdot \Delta\hat{\beta})) = 0 \text{ for some } m \in \mathcal{G}_k \text{ with } k \notin \mathcal{A}_\lambda\},$$

where we take the infimum over an empty set to be $+\infty$.

(e) Set $t = t + 1$, $\delta = \min\{\delta_A, \delta_I, \delta_R, \delta_U, \delta_S, \lambda_t\}$, $\lambda_{t+1} = \|c_k(\hat{\beta}(\lambda_t) + \delta \cdot \Delta\hat{\beta})\|_\infty$ and $\hat{\beta}(\lambda_{t+1}) = \hat{\beta}(\lambda_t) + \delta \cdot \Delta\hat{\beta}$.

**Acknowledgments.** The authors would like to gratefully acknowledge the comments by Nate Coehlo, Jing Lei, Nicolai Meinshausen, David Purdy and Vince Vu.

DEPARTMENT OF STATISTICS
UNIVERSITY OF CALIFORNIA, BERKELEY
367 EVANS HALL
BERKELEY, CALIFORNIA 94720
USA
E-MAIL: pengzhao@stat.berkeley.edu
gvrocha@stat.berkeley.edu
binyu@stat.berkeley.edu